\documentclass[a4paper,11pt]{article}
\usepackage{amsfonts,amssymb,amsthm,amsmath,fancyhdr}

\newcommand{\annsection}[1]{\stepcounter{section}  \noindent {\bf \thesection. #1. }}
\newcommand{\annsubsection}[1]{\stepcounter{subsection} \noindent{\it \thesubsection. #1. }}

\title{On Local Poincar\'e via Transportation}
\author{Max-K. von Renesse\footnote{Currently visiting  Courant Institute with the support of the Alexander von Humboldt Foundation (AvH).}\\
Technische Universit\"at Berlin\\
{\tt mrenesse\@@math.tu-berlin.de}\\
  }

\begin{document}

 \maketitle \abstract{It is shown that curvature-dimension bounds $CD(N,k)$ 
for a metric measure space $(X,d,m)$ in the sense of Sturm imply a weak $L^1$-Poincar\'e-inequality provided $(X,d)$ has $m$-almost surely no branching points.} \\


\def\N{\mathbb N}
\def\M{\mathbb M}
\def\fdp{f\circ d_{p}}
\def\R{\mathbb R}
\def\DE{D(\mathcal E)}
\def\DEC{\mbox{D}_c(\mathcal E)}
\def\mint{\,\mathop{\makebox[0pt][c]{$\int$}\makebox[0pt][c]{--}\,}}
\def\mintl{\!\!\mathop{\makebox[0pt][c]{$\int$}\makebox[0pt][c]{--}\!}\limits}
\def\suml{\sum\limits}
\def\Er{\mathcal E^r}
\def\E{\mathcal E}
\def\Erw{\mathbb E}
\def\intl{\int\limits}
\newcommand{\Br}[1]{B_r(#1)}
\newcommand{\Bbr}[1]{\mathbb B_r(#1)}
\newcommand{\bnk}[1]{b_{n,k}(#1)}
\newcommand{\bn}[1]{b_{n}(#1)}
\def\phx{\phi_x}
\def\Mkx{\M_k(x)}
\def\bm{{\mathop{\mbox{\normalfont \,vol}}}}
\def\dbm{{{\mbox{\normalfont \tiny d \!\!vol}}}}
\def\dbmn{{{\mbox{\normalfont d \!\!vol}}}}
\def\Qrx{q_r(x)}
\def\Regset{X\setminus S_X}
\def\volsing{\mathbb S_X}
\def\angle{\sphericalangle}
\def\Eomega{\E_\Omega}
\def\Deomega{\mathcal D(\Eomega)}
\def\Meas{\cal M}
\def\1{\,{\makebox[0pt][c]{\normalfont    1}
\makebox[2.5pt][c]{\raisebox{3.5pt}{\tiny {$\|$}}}
\makebox[-2.5pt][c]{\raisebox{1.7pt}{\tiny {$\|$}}}
\makebox[2.5pt][c]{} }}
\def\one{\1 }
\def\one{\1 }
\def\B{\mathcal B}
\newcommand{\infl}[1]{\inf\limits_{#1}}
\newcommand{\gsk}[1]{\left\{ #1 \right\}}
\newcommand{\norm}[1]{\left\| #1 \right\|}
\def\P{{ /\!\!/}}
\def\ul{\underline}
\def\H{\mathcal H}
\def\for{\mbox{ for }}
\def\Ricc{\mathop{\mbox{\normalfont Ricc}}}
\def\Sec{\mathop{\mbox{Sec}}}
\def\ra{\rangle}
\def\la{\langle}

\def\sth{\,|\,}
\def\Hess{\mathop{\mbox{{\normalfont Hess}}}}
\def\einh{\frac{1}{2}}
\def\ol{\overline}
\def\abs{~\\~\\}
\def\sabs{~ \smallskip ~}
\def\k{,\,}
\def\supp{\mathop{\mbox{supp}}}
\def\Curv{\mbox{Curv}}

\def\cut{\mbox{\normalfont Cut}}
\def\mt{\rightarrow}
\newcommand{\lgauss}[1]{{\lfloor #1 \rfloor}}
\def\falls{\mbox{ if }}
\def\sonst{\mbox{ else}}
\def\grad{\mbox{grad}}
\def\div{{{\mathop{\,{\rm div}}}}}
\def\diam{{{\mathop{\,{\rm diam}}}}}
\def\Id{{\mathop{{\bf 1}_{\small \R^d}}}}
\def\k{,\,}
\def\dist{{\mbox{dist}}}
\def\Ind{{\mathop{\mbox{I}}}}
\def\dil{{\mathop{\mbox{dil}}}}
\def\H{\mathbb H}
\def\mtm{M \times M}
\def\BV{\mbox{\normalfont BV}}
\def\gij{{g_{ij}}}
\def\gije{{g_{ij}^\epsilon}}
\def\Lip{\mathop{\mbox{Lip}}}
\def\loc{{\mbox{\scriptsize loc}}}
\renewcommand{\bullet}{{\mathbf \cdot  }}
\def\Cut{\mbox{\normalfont Cut}}
\def\dist{\mbox{\normalfont dist}}
\def\sign{\mbox{\normalfont sign}}
\def\pr{\mbox{\normalfont pr}}
\def\Z{{\mathbb Z}}

\def\ioe{\frac{1}{\epsilon}}
\def\fe{f_{\epsilon}}
\def\Ebr{\E^{b,r}}
\def\Dco{D_c(\E_\Omega)}

\newenvironment{bew}{\begin{proof}}
{\end{proof}}


%
%
\newtheorem{thm}{Theorem}
\newtheorem{rthm}{Theorem}
\newtheorem{satz}{Proposition}
\newtheorem{lem}[]{Lemma}
\newtheorem{cor}[]{Corollary}
\newtheorem{defn}[]{Definition}

\def\P{{\mathcal P}}
\def\Id{{\mathop{\rm Id}}}
\parindent0cm
 
\annsection{Introduction} From analysis on manifolds and  metric measure spaces $(X,d,m)$ the fundamental importance of Poincar\'e-type inequalities  for the regularity of harmonic, Lipschitz or Sobolev functions is known (cf. 
\cite{MR97k:31010, MR2000g:53043, MR1683160} and \cite{MR1800917,MR1872526,MR2039660}). In this note we show that metric measure spaces $(X,d,m)$  with  upper dimension-lower Ricci curvature $CD(N,k)$ bounds in the generalized sense of Sturm \cite{MR2123035}  (cf. \cite{LV04} for $CD(N,0)$) support a weak local $L^1$-Poincar\'e inequality provided the choice of geodesics by which mass travels between cut points in $(X,d)$ can be made symmetric. By some additional argument on the support of the measure $m$ we show that this will be the case if we assume moreover that the underlying metric admits $m$-almost surely no branching points. 

\smallskip

The main step is to establish a version of the {\it segment inequality} of Cheeger-Colding \cite{MR1405949}
 which in particular implies the Poincar\'e inequality in the sense of upper gradients.
   What is needed for the proof is only a quantitative version of the Bishop-Gromov volume comparison 
theorem which is called $(N,k)$-measure contraction property in \cite{ohta05}
and which is implied by the Lott-Villani-Sturm  dimension-curvature bounds. 

\smallskip

The assumption on the $m$-almost sure absence of branching points of $(X,d)$ may be quite restrictive.
However, as a first step in understanding the full meaning of curvature-dimension bounds for metric measure spaces it may be a useful task to study the regularity of admissible spaces $(X,d,m)$ witout branching points and study their relation with Alexandrov spaces, for instance. 

\smallskip

\annsection{Preliminary study: Proof  of the segment inequality on smooth manifolds by  mass transportation}

\smallskip

For illustration let us derive the segment inequality in the smooth case using the language of mass transportation. This approach proved  recently very useful for the generalization of certain concepts in smooth Riemannian analysis to general metric measure spaces \cite{VRS, MR2123035, LV04, Sturm05}.\\

We review some standard terminology first. Throughout this note we call a curve $\gamma:[0,1]\mt (X,d)$ in a metric space a {\it geodesic} (segment) if $d(\gamma(s),\gamma(t)) = d(\gamma(0),\gamma(1))|s-t|$ for all $s,t\in [0,1]$. For $A,B \subset   X$ we define the set $\Gamma(A,B)$ as the collection of all geodesics $\gamma$ 
with $\gamma(0)\in A$ and $\gamma(1)\in B$. For $x,y\in X$ any $\gamma \in \Gamma(x,y)$ will be denoted $\gamma_{xy}$ which may not be unique. If $X= M$ is a Riemannian manifold $(M^n, g)$ then $d$ is the intrinsic metric induced by $g$. \\

{\bf Propostion. (\cite{MR1405949})} {\it Let $(M^n,g)$ be a smooth Riemannian manifold with Ricci curvature $\mathop{\rm Ricc}_M \geq (n-1) k$, $k\in \R$. Let $A_1, A_2 \subset B_R$ be measurable subsets contained in  a geodesic $R$-Ball and let $g:B_{2R}\mt \R$ be a nonnegative measurable function, then 
\[ 
\int_{A_1} \int_{A_2}
\int_0^{  1}  g(\gamma_{xy}(s)) dt  dx dy 
\leq C_\kappa(n,D) (|A_1|+ |A_2|)\int_{B_{2R}} g(z) dz, \]
with
\[ 
 C_\kappa(n,D) = \sup_{  t \in [\frac 1 2 , 1]}
 t \left( \frac { s_k( k D)}{ s_k( t D )}\right)^{n-1}
\] 
and $D=D(A_1,A_2)= \sup_{(x,y)\in A_1\times A_2} d(x,y) \leq {\rm Diam}(A_1\cup A_2)$.
}\\

Here and in the sequel $s_k$ denotes  the usual Sturm-Liouville function 
\[  s_k(t) = \left\{
\begin{array}{ll}
\sin(\sqrt k t) & \falls k >0 \\
t & \falls k =0\\
\sinh(\sqrt{ -k }t) & \falls k <0. \\
\end{array}
\right.
\]

Let us recall the necessary basic facts from optimal mass transportation theory we need 
 (cf. \cite{MR1865396} and  \cite{MR1964483} as general reference).
 
 \smallskip
 
Let $\mu_0$ and $\mu_1$ be two probability measures on a Riemannian manifold $M^n$ with $\Ricc(M) \geq (n-1)\kappa$ and let $\tau_t : M\mt M$, $t \in [0,1]$ the optimal transportation map associated to the $L^2$-Wasserstein metric with the squared distance as cost function. General theory says that for $\mu_0 \ll d{\rm vol}_{g}$ this map is of the form 
\[ \tau_t (x) = \exp_x(-t\nabla\Psi(x))\]
where $\psi: M\mt \R$ is a $d^2/$-concave function (i.e $\psi$ it is the inf-convolution of another function $\phi$ with respect to the potential $d^2/2$).\\

The main ingredient of the proof given below is the concavity estimate for the n-th root of the Jacobian of $\tau_t$ with respect to $t$ which is one of the main results in the fundamental paper \cite{MR1865396}.
 It reads as follows. If $J_t(x):= \det d\tau_t(x)$ denotes the Jacobian determinant of the map $\tau_t $ in $x$ then  
\[ J_t ^{1/n}(x) \geq \tau_{k}^{(1-t)}(x)J_0^{1/n}(x) +  \tau_{k}^{(1-t)}(x)J_1^{1/n}(x)\]
with   \[
 \tau_{\kappa}^{(t)}(x):= t^{1/n}\left(\frac{s_k( t d(x,\tau_1(x))}{s_k(  d(x,\tau_1(x))}\right)^{1-1/n}=:\tau_{k,n}(t(d(x,\tau_1(x))).\]
  
Here $d(\tau_1(x),x)$ is the distance between a point $x$ and its target point $ \tau_1(x)$. 
It is obtained using the structure of $\tau_t$ and Jacobi-field estimates together with the central arithmetic-geometric mean inequality on nonnegative matrices (cf. \cite{Sturm05} for a nice presentation).  \\

{\it Proof of the Riemannian segment inequality.} Let  be $A_1 $ and $A_2$ be two sets which are both embedded in a larger ball $B_R$ and let 
 $y \in A_2$ by a point. Let $\mu_t={\tau_t}_*\mu_0$ the Wasserstein geodesic connecting   $\mu_0= \frac 1 {|A_1|}dx_{|A_1} $ on $A_1$ with $\delta_y$, the Dirac measure in $y$. From the structure of  it is clear that each $x \in A_1$ travels along a geodesic connecting $y $ with $x$ (i.e. $\psi(x)= d^2_y(x)/2$ in $\tau_t$).  Since the cut locus of $y$ has measure zero we can assume below that there is only one such geodesic for each pair $(x,y) \in M\times M$. Let $g: B_{2R}\mt \R$ be nonnegative, then 
\begin{align*}
\frac 1 {|A_1|} & \int_{A_1} \int_0^{\frac 1 2}  g(\gamma_{xy}(s)) dt  dx  =\frac 1 {|A_1|} \int_{A_1} \int_0^{\frac 1 2} g(\tau_t^y(x)) dt  dx \intertext{which by the general integral transformation formula equals}
& =\int_0^{\frac 1 2} \int_{A_{1_t}}  g(z) \mu_t (z) dz dt = \int_0^{\frac 1 2} \int_{A_{1_t}}  g(z) \frac{m_0 (x(t,z))} {\det d\tau_t(x(t,z)) } dz dt
\end{align*}
 {with $A_{1_t}:= \{\gamma_{yx}(t)\,|\, y \in A_1\}$ 
and where we used the Jacobi identity for $m_t = \frac {d\mu_t}{dx}$
\[  m_t(\tau_t(.)) \det d\tau_t (.) = m_0(.)\]
with  $x(t,z)=\tau_t^{-1}(z)$ being the origin of the transport ray which hits $z$ at time $t$. By the 
 Jacobian concavity 
\[ \frac 1 {(\det d\tau_t(x(t,z))}   \leq   \left( \tau_{\kappa}^{(1-t)}(x(t,z))   
 + \tau_{\kappa}^{(t)}(x(t,z))    \det d\tau_1^{1/n}(x(t,z)) \right)^{- n} \] 

where in the present case  $\det d\tau_1(x(t,z))=0$ since $\mu_1 = \delta_y$. (To see this  approximate $\mu_1=\delta_y$ by the family  $\mu^\epsilon_1 = \frac 1{|B_\epsilon(y)|} dx _{|B_\epsilon (y)}$ and use   Jacobian identity for the density $m^\epsilon_1$ of  $\mu_1^\epsilon$
\[  m^\epsilon_1(z)= \frac 1{|B_\epsilon(y)|} \one_{B_\epsilon(y)}(z) = \frac {m_0(x(z))}{\det d\tau_1 ^\epsilon(x(z))}= \frac {1/ {|A_1|}\one_{A_1}(x(z))}
{ \det d\tau^\epsilon_1(x(z))}\]
which implies   $ \det {\det d\tau_1 ^\epsilon(x(z))} \equiv \frac {|B_\epsilon(y)|}{|A_1|}$ for all $z \in B_\epsilon(y)$, 
equivalently  $  \det d\tau_1 ^\epsilon(x) \equiv \frac {|B_\epsilon(y)|}{|A_1|}$ for all $x \in A_1$ or $ \det d\tau_1 ^\epsilon(x(t,z)) \equiv \frac {|B_\epsilon(y)|}{|A_1|}$ for all $z \in A_{1_t}$, and let $\epsilon \to 0.)$ 
Since ${m_0 (x(t,z))}= \one_{A_{1_t}}(z) /|A_1|$ we arrive at the estimate 
\begin{align*}
\frac 1 {|A_1|}  \int_{A_1} \int_0^{\frac 1 2} & g(\gamma_{xy}(s)) dt  dx \\
& \leq 
  \int_0^{\frac 1 2}  \int_{A_{1_t}} g(z){m_0 (x(t,z))} \left( \tau_{\kappa}^{(1-t)}(x(t,z))\right)^{- n}dzdt \\
     & \leq \frac 1 {2 |A_1|} \sup_{
 \substack{  t \in [\frac 1 2 , 1] \\x \in A_1}
 }
     \left( \tau_{\kappa}^{t}(x)   \right)^{- n} 
     \int_{B_{2R}} g(z)dz  , 
\end{align*}
where we used  $A_1, A_2 \subset B_R$.
Integration  respect to $y \in A_2$ yields
\[ 
\int_{A_1} \int_{A_2}
\int_0^{\frac 1 2}  g(\gamma_{xy}(s)) dt  dx dy 
\leq C_\kappa(n,D) \frac{|A_2|}2 \int_{B_{2R}} g(z) dz \]
where $
 C_\kappa(n,D)  $ is the constant as defined in the statement of the proposition. 
 This is because  the  expression  for $  \tau_{\kappa}^{t}(x)$ is monotone increasing in $d(x,y)$. - Using the symmetry of the integral estimate we can bound the expression 
 \[ 
\int_{A_1} \int_{A_2}
\int_{\frac 1 2}^{1}  g(\gamma_{xy}(s)) dt  dx dy 
\leq C_\kappa(n,D) \frac{|A_2|}2 \int_{B_{2R}} g(z) dz \] 
by repeating the previous arguments to the corresponding integral over the time interval $[0,{\frac 1 2}]$ when $A_1$ and $A_2$ are interchanged (see also section 3.) which by adding the two estimates 
concludes the proof. \hfill$\Box$

\smallskip

In the estimate above we defined the geodesic to be parameterized on the unit interval. Using unit speed  parameterization  it reads  \[ 
\int_{A_1}\int_{A_2} \int_0^{ {d(x,y)} }  \frac{g(\gamma_{xy}(t))}  {d(x,y)}  dt dx dy
\leq \frac 1 2 C_\kappa(n,D)(|A_1|+|A_2|) \int_{B_{2R}} g(z) dz.
\]
~

\smallskip

 \annsection{Segment Inequality  on metric measure spaces with transportation lower Ricci bounds}\\

\annsubsection{Measure contraction}

\smallskip

We start out from the following definition which is a quantified version of the measure contraction property formulated in \cite{Sturm98}.\\

{\bf Definition. (\cite{ohta05})} {\it A metric measure space $(X,d,m)$ 
is having the
{ $(N,k)$ measure contraction property} if
for each pair $(x,M) \in X\times   2^X$ there exists a probability 
measure $\Pi$ on the 
set of geodesics $\Gamma(x,M)=\{ \gamma_{xy}|y \in A\}$ with  
$e_{0*}\Pi =  m_M:= \frac 1{|M|} m_{|M}$ and $e_{1*}\Pi =\delta _x$ such that
\[  dm \geq e_{t*} \left( (1-t) \left\{\frac {s_k((1-t) \ell (\gamma)}{s_k(\ell(\gamma))}\right\}^{N-1}\mu(M) d\Pi\right),\]
where $e_t: \Gamma \mt X$ with $e_t\gamma := \gamma(t)$, $t\in [0,1]$ is the evaluation map. 
}~\\

Let us abbreviate this property by $MCP(N,k)$. Its meaning  is the following. 
Disintegrating the measure $\Pi$ with respect to the evaluation map $e_0$ and using the condition that $e_{0*}\Pi= m_M$ we obtain a mixing representation  $\Pi(d\gamma) = \lambda_y(d\gamma) m_M(y)$ where the measures  $\lambda_y$ are  supported on $\Gamma(y,x)$.
Moreover,  the measures $\lambda_{yx}$ are determined uniquely by $\Pi$ for $m$-almost $y$   and vice versa.

\smallskip

Let now $M_t:= e_t \Gamma(M,x)$ be the set hit by all geodesics from $M$ to $x$ then
  the statement above is equivalent to 
\begin{equation}
  m(Z) \geq (1-t)  \intl_{M} \intl_{\Gamma(y,x)} \one_Z(\gamma_t)  \left\{\frac {s_k((1-t) d(x,y))}{s_k(d(x,y)}\right\}^{N-1} \lambda_y(d\gamma) m(dy)\label{mcp}\tag{MCP}
 \end{equation}
 for all measurable $Z$ with w.l.o.g $Z \subset M_t$, since for $Z \subset X\setminus M_t$ the right hand side is zero. - Written in this form the $MCP$-condition gives a lower bound for the concentration of $m$ under the generalized homothetic map defined  by the Markov kernels $\Lambda^y_{t}(dx) = e_{t*}(\lambda_{y})(dx)$ (cf. next section). It may also be seen as a requirement on the minimal 'mean spreading' of  all geodesics to $x$, where the mean is taken with respect to the collection of weights $(\lambda_y)$ and $m$.   \\

The relevance of the $(N,k)$-measure contraction property is its robustness with respect to the measured Gromov-Hausdorff convergence (cf. \cite{ohta05}). Moreover, it is implied by the generalized dimension-curvature bounds defined recently by Sturm.

\smallskip

{\bf Definition. (\cite{Sturm05})} {\it A metric measure space $(X,d,m)$ satisfies the curvature dimension condition $CD(N,k)$ if for each pair $\nu_0 , \nu_1 \in \mathcal P_2(M,d,m)$ there exists an optimal $d^2$-coupling $q \in \mathcal P(X\times X)$ and a  geodesic $\Gamma:[0,1]\mt  P_2(X,d,m)$ connecting $\nu_0$ and $\nu_1$ such that 
\begin{align*}
S_{N'}(\Gamma(t)|m) \leq - \int_{X\times X}& \left[ 
\tau_{K,N'}^{(1-t)}(d(x_0,x_1)\cdot \rho_0^{-1 /N'}(x_0) \right.\\
 & ~~~~~~\left.  +\tau_{K,N'}^{t}(d(x_0,x_1)\cdot \rho_1^{-1 /N'}(x_1)\right] dq(x_0,x_1) 
 \end{align*}
where $\rho_i = d\mu_i/dm$ are the respective densities of $\mu_i$, $S_{N'}(\mu|m)= - \int_X \rho^{-1/N'} dm $ is the Renyi entropy of a measure $d\mu = \rho dm $ with respect to $m$ and the functions $\tau_{K,N'}$ are defined as above.} 

\smallskip

The following is easily verified (see \cite{Sturm05}, proposition 1.7.iv.]).

\smallskip
 
{\bf Proposition.} {\it   $CD(N,k)$ implies $MCP(N,k)$  for a metric measure space $(X,d,m)$.}

\smallskip

From the results in \cite{MR1865396} one deduces   that smooth Riemannian manifolds $(M^n,g)$ with Ricc$(M^n,g)\geq (n-1)k$ satisfy the $CD(n,k)$ condition. - One of the main results in Sturm's theory  is 
the stability of the $CD(N,k)$ bounds with respect to convergence.

\smallskip
 
{\bf Theorem. (\cite{Sturm05})} {\it For fixed $N,k \in \R$ the set of metric measure spaces satisfying the $CD(N,k)$ is closed with respect to measured Gromov-Hausdoff convergence.}

\smallskip

In particular, measured Gromov-Hausdorff limits of $CD(N,k)$-spaces  will satisfy the $(N,k)$-measure contraction property. (The corresponding  stability result of the $(N,K)$-MCP condition alone is a little stronger and is obtained in \cite{ohta05}  based on   \cite{LV04}).\\

\annsubsection{Segment inequality for $(X,d,m)$}

\smallskip

For the precise formulation of the subsequent results we need a little more notation.
Let  $B \subset X$ be a set and   $t\in [0,1]$. Then we define  the following set valued 
{\it geodesic contraction map} in direction $B$
\[  \Gamma_t(.,B): 2^X \mt 2 ^X ; \quad \Gamma_t(A,B) := e_t (\Gamma(A;B)).\] 
 By abuse of notation for $t\in [0,1]$  we define also 
$ 
\Phi_{-t}(A,B)=\{ z \in X \,| \, \exists \, b\in B, \gamma_{zb}: \gamma_{zb}(t) \in A\}$ the inverse  $A$ 
with respect to geodesic contraction in direction $B$. When $B=\{b\}$ we write $A_t(b):= \Gamma_t(A,\{b\})$ and $A_{t}^{-1}(b):= \Gamma_{-t}(A,\{b\})$. \\

Note that in the $(N,k)$-MCP statement above the transference plan $\Pi$ and thus the measures $(\lambda_{yx})_{y\in M}=(\lambda^M_{yx})_{y\in M}$ may depend on $M$.  Let us say that the family of measures $(\lambda^M_{xy})_{x,y \in X}$ is {\it symmetric} in $(x,y)$ if 
\[ 
\lambda^M_{xy} (d\gamma) = \lambda^M_{yx}(d\ol\gamma),\] 
where $\gamma \mt \ol \gamma$ is the inversion map $\ol \gamma(t)=\gamma(1-t)$, $t \in [0,1]$.\\

{\bf Proposition.} {\it Let $(X,d,m)$ satisfy the $(N,k)$-measure contraction property and assume that 
the map $(x,y) \mt \lambda^M_{xy} \in \P(\Gamma_{xy})$ can be chosen  
to be symmetric for $m\times m$ almost every pair $(x,y)$  and for 
some set $M\subset X$. Then  for two measurable subset $A_1, A_2 \subset B_R$ contained in a geodesic ball $B_R\subset M$ and 
$g: B_{2R}\mt \R$ nonnegative 
\begin{align*}
\intl_{A_1}\intl_{A_2} \intl_{\Gamma(x,y)}\intl_0^{ 1} & {g(\gamma_{xy}(t))}   dt   \lambda^M_{xy}(d\gamma)
 m(dx)m(dy)\\
 & 
\leq \frac 1 2 C_\kappa(n,D)(|A_1|+|A_2|) \intl_{B_{2R}} g(z) m(dz).
\end{align*}
}

\smallskip

{\it Proof.} We write the  $(N,k)$-measure contraction 
inequality relative to the ambient set $M$ 
yet in another form, namely for all $(x,t) \in X\times [0,1]$

\[ m \geq (\tau_t^x\cdot \Lambda^x_t) * m, 
\]
where the sign $*$ means convolution  with the transition kernel 

\[  (\tau_t^x\cdot \Lambda^x_t) (y,dz) = \tau(x,y) \Lambda_t^x (y,dz) \]
with the symmetric function 
$\tau(x,y)= t \left\{\frac {s_k((1-t) d(x,y))}{s_k(d(x,y)}\right\}^{N-1}$  
and the Markov kernel 
\[ 
 (y,Z) \mt  \Lambda_t^x (y,Z) = \int_{\Gamma(y,x)} \one_{Z}(\gamma_t)\lambda_{yx}(d\gamma),  \quad y \in X, Z \subset X.\]
Here we omitted the upper index  for $ \lambda_{yx}= \lambda^M_{yx}$ as we shall in the rest of the proof because $M$ is fixed. 
\\
 
With this notation the $(N,k)$-measure contraction 
inequality is written 
\[\int_{Z} g(z) m(dz) \geq \int_{X}  \int_ {Z} 
g(y) \tau_t (x,z) \Lambda^x_t(z,dy) m(dz)\]
for all measurable $Z \subset X$ and nonnegative measurable $f:X\mt \R$. Note that it suffices to take the outer integral on the set $ 
Z_t^{-1}(x)$. 
Since $ d(z, x) =  {d(\gamma_{zx}(t),x)}/{(1-t)}$ for all $z\in X$  
we have 
\[ 
\inf_{z  \in  Z_t^{-1}(x)} \left[ \tau_t (x,z)\right] = \inf_{z  \in  Z}\left[ \tau_t \bigl({d(x,z)}/({1-t})\bigr)\right]
\]
such that the estimate 
\[\int_{Z} g(z) m(dz) \geq  \inf_{z  \in  Z}\left[ \tau_t (d(x,z)/(1-t))
\right] \int_{X}  \int_ {Z} 
g(y)\Lambda^x_t(z,dy) m(dz).\]
is obtained. For $A  \subset  X$ set $Z = A_t(x)$, then  
under the assumption $A\subset B_{R}, x\in B_R$ it follows $A_t(x) \subset B_{2R}$, thus
\[\int_{B_{2R}} g(z) m(dz) \geq  \inf_{z  \in  A}\left[ \tau_t (z,x)
\right] \int_{X}  \int_ {A_t(x)} 
g(y)\Lambda^x_t(z,dy) m(dz).\]

Integration of this inequality with respect to time yields 
%
\[ \int_{0}^{\frac 1 2}
 \int_{X} \int_ {A_t(x)} g(y)  \Lambda^x_t(z,dy) m(dz)dt \leq \frac 1 2  \sup_{\substack{ t\in [ 0,\frac 1 2  ]\\ z \in {A_1}}}[ \tau_t^{-1} (x,z) ]
  \int_ {A_1}   g(z) m(dz). \] 

Since 
\begin{align*}
 \intl_X\intl_{A_t(x)} g(y)  \Lambda^x_t(z,dy)m(dz) & =  \intl_{\Gamma_{-t}(A_t(x),x)}\intl_{A_t(x)} g(y)  \Lambda^x_t(z,dy)m(dz) \\
& \geq  \int_A \int_{\Gamma(X,X)} g(\gamma(t))\lambda_{zx}(d\gamma)m(dy)
\end{align*}

the estimate 
\[ 
  \int_ {A} \int_{\frac 1 2}^1\int_{\Gamma(X,X)} g(\gamma(t))\lambda_{zx}(d\gamma)dt m(dz)  \leq \frac 1 2  \sup_{\substack{ t\in [\frac 1 2 , 1 ]\\ z \in {A}}}[ \tau_t^{-1} (x,z) ]
  \int_ {A} g(z) m(dz). \]  
 is obtained. Now put $A=A_1$ and integrate the last inequality   with respect to $x \in A_2$. This leads to 
\begin{align*}
   \int_ {A_2}  \int_ {A_1} 
   \int_{\frac 1 2}^1\int_{\Gamma(X,X)} & g(\gamma(t))\lambda_{zx}(d\gamma)dt m(dz) m(dx)\\
    & \leq \frac 1 2  |A_2| 
   \sup_{\substack{ t\in [ 0,\frac 1 2  ]\\ (z,x) \in A_1\times A_2 }}[ \tau_t^{-1} (x,z) ]
  \int_ {A_1}   g(z) m(dz) \\
&  \leq \frac 1 2  |A_2| C_k(N,D)
     \int_ {A_1\cup A_2}   g(z) m(dz),
  \end{align*}
with $C_k(N,D)$ and $D= D(A_1,A_2)$ defined as in the smooth case above.\\

Finally, interchanging the roles of $A_1$ and $A_2$ and using the symmetry of the measures $(\lambda_{xy})$ we obtain
\begin{align*}
   \frac 1 2  |A_2| C_k&(N,D)
     \int_ {B_{2R}}   g(z) m(dz)\\
     &  \geq  \int_ {A_1}  \int_ {A_2} 
   \int_0^{\frac 1 2}\int_{\Gamma(X,X)} g(\gamma(t))\lambda_{zx}(d\gamma)dt m(dz)m(dx)  \\
       &  =  \int_ {A_1} 
\int_ {A_2}  \int_0^{\frac 1 2}\int_{\Gamma(X,X)} g(\gamma(t))\lambda_{xz}(d\ol \gamma)dt m(dx)m(dz) \\ 
 &  = \int_ {A_1} \int_ {A_2}   
   \int_0^{\frac 1 2}\int_{\Gamma(X,X)} g(\ol \gamma(1-t))\lambda_{xz}(d\ol \gamma)dt m(dx)m(dz) \\
 &  = \int_ {A_1} \int_ {A_2} 
   \int_{\frac 1 2}^1\int_{\Gamma(X,X)} g(\ol \gamma(t))\lambda_{xz}(d\ol \gamma)dt m(dx)m(dz) \\
 & = \int_ {A_1} \int_ {A_2}  
   \int_{\frac 1 2}^1\int_{\Gamma(X,X)} g( \gamma(t))\lambda_{xz}(d \gamma)dt m(dx)m(dz)\\
   & = \int_ {A_2}  \int_ {A_1} 
  \int_{\frac 1 2}^1\int_{\Gamma(X,X)} g( \gamma(t))\lambda_{zx}(d \gamma)dt m(dz)m(dx). 
\end{align*}

Adding this inequality to the first one the claim is established. \hfill $\Box$\\

{\bf Corollary.} {\it  The assertion of the proposition above is true in particular when the cut-locus $C_x:=\{ y \in X \,|\, \#\Gamma(x,y)  \geq 2 \} $ satisfies $m(C_x)=0$ for $m$-a.e. $x \in X$.}  

\smallskip

{\it Proof.} If $x \not \in C_y$ then $\lambda_{xy}=\delta_{\gamma_{xy}}$ for the unique $\gamma_{xy} \in \Gamma(x,y)$. This forces $\lambda_{xy}$ to be symmetric $m\times m$-almost surely. \hfill $\Box$\\

For the following version of the Poincar\'e inequality recall that for  $f: X\mt \R$ 
the function $g: X\mt \R_+$ is called an {\it upper gradient} if
\[  |f(x)-f(x)| \leq \int_0^{d(x,y)} h(\gamma_s) ds \]
for any unit speed geodesic connecting $x$ and $y$.  

\smallskip

{\bf Corollary. ($L^1$-Poincar\'e-inequality)} {\it Under the conditions above let $h$ be an upper gradient of $f$ then 
\[
\int_{B_R}\int_{B_R} \frac{|f(x)-f(y)|}{d(x,y)} m(dx)m(dy)
\leq |B_R| C_{n,k}(D) \int_{B_{2R}} h(x)m(dx).
 \]}

\smallskip
 
{\it Proof.}  In order to prove this inequality  for each pair $(x,y)$ let $\lambda_{xy}$ be the associated 
measure 
from proposition above, then the assertion follows from   
\[ \frac{|f(x)-f(y)|}{d(x,y)} =  \int_{\Gamma_{xy}} |f(\gamma(0))-f(\gamma(1))| \lambda_{xy}(d\gamma)
\leq   \int_{\Gamma_{xy}} \int_0^1 h(\gamma_s) ds \lambda_{xy}(d\gamma)\] 
which can be inserted into the segment inequality. \hfill $\Box$\\

{\it Examples.} Consider the Banach space $(X,d)= (\R^n , \norm{.}_p)$, $p \in (1,\infty]$ equipped with $m= \lambda^n$ the $n$-dimensional Lebesgue measure, where $\norm{x}_p = (\sum_{i=1}^n  |x_i|^p )^{1/p}$. The geodesics are straight Euclidean line segments $\gamma_{xy}= x+ t(x-y)$, hence $C_x=\emptyset $ for all $x \in X$. Obviously $(\R^n , \norm{.}_p,\lambda^n)$ satisfies $MCP(0,n)$. - The case $(\R^n , \norm{.}_1,\lambda^n)$ is a little more interesting, since here $C_x=X$ for all $x \in X$.  However, choosing $\lambda_{xy}(d\gamma) =\delta_{(x+ t(x-y))}(d\gamma)$ for all $x,y \in X$ the  $MCP(0,n)$ property remains true. Since this choice of $\lambda$ is symmetric, the segment and Poincar\'e inequalities hold.\\

\annsubsection{Extendable geodesics and branching points}\\
  
  {\bf Definition.} {\it Let $(X,d)$ be a metric set and $x \in X$. Define $I_p:= \{ \gamma_{xp}(t)| t\in (0,1] , x\in X\}$ as all points $x$ which are connected to $p$ by at least one extendable geodesic segment and let $T_p=X\setminus I_p$.}\\

{\bf Proposition.} {\it Let $(X,d,m)$ satisfy an $(N,k)$-measure contraction property then
$m(T_p)= 0$ for all $p\in X$.}

\smallskip

{\it Proof.}  
The proof is an adaption of the idea behind 
proposition  3.1. in \cite{MR1274133}. - Note first that $I_p = \bigcup _{t\in (0,1]} X_{t}(p)$ where the sets $X_{t}(p) $ are monotone decreasing  for $t\in [0,1]$, i.e. $X_{t}(p)\subset X_{s}(p)\subset $ for $s\leq t$.
 Let $ A \subset X $ be a measurable bounded set and 
choose the  weight functions $\lambda_{xy}= \lambda^M_{xy}$ for $M$ large enough such that 
$A \subset M_{-t}(p)$ for all $t\in [0,1]$.
 The measure contraction inequality at time $t $ 
applied to the set $A\cap X_{t}(p) $ yields 
\begin{align*}
m(  A\cap X_{t}(p) ) 
& \geq  (1-t) \int_X  \left(\frac {s_k((1-t) d(p,y))}{s_k(d(p,y))}\right)^{N-1} \lambda_{yp}^t( A\cap X_{t}(p) ) m(dy)\\
& = (1-t) \int_X  \left(\frac {s_k((1-t) d(p,y))}{s_k(d(p,y))}\right)^{N-1} \lambda_{yp}^t( A ) m(dy)
\end{align*}
 {because $\lambda_{yp}^t( A\setminus X_{t}(p) )=0 $ for all $y\in X$. Hence for  $s \leq t$ we obtain from the monotonicity of the family $(M_s(p))_{s\in [0,1]}$ and } 
\begin{align*}
m(  A\cap X_{s}(p) ) 
& \geq  (1-t) \int_X   \left(\frac {s_k((1-t) d(p,y))}{s_k(d(p,y))}\right)^{N-1} \lambda_{yp}^t( A  ) m(dy).
\end{align*}
 {Sending $s \to 0$, since  $X_s(p) \nearrow I_p$ for $s\searrow0$ by monotone convergence}
\begin{align*}
m(  A\cap I_p ) 
& \geq  (1-t) \int_X  \left(\frac {s_k((1-t) d(p,y))}{s_k(d(p,y))}\right)^{N-1} \lambda_{yp}^t( A  ) m(dy).
\end{align*} 
Finally, upon  sending $t\to 0$ and using $\lambda^0_{yp}(dz)= \delta_y(dz)$ 
 dominated convergence yields
\[m(  A\cap I_p ) 
 \geq  m(A),\]
which by the arbitrariness of $A$ establishes the claim that $I_p$ has full $m$-measure.   \hfill ~$\Box$\\

A geodesic tripod is a space $(X,d)$ obtained by metric 
gluing of three segments in one common endpoint which is called the center.

\smallskip  
  
{\bf Definition.} {\it  
A point $p\in X$ in a metric space 
$(X,d)$ is called  branching point if it is the center an embedded geodesic tripod. } 

\smallskip
 
{\bf Corollary.} {\it If $(X,d)$ admits   no branching points $m$-almost surely and $(X,d,m)$ satisfies an $MCP(n,k)$-property then 
$m(C_x)=0$ for  $x \in X$. Also, the segment and Poincar\'e inequalities hold for $(X,d,m)$ in this case.} 

\smallskip

{\it Proof.} Since for $m$-almost all $y  \in C_x$ at least one $\gamma_{xy} \in \Gamma(x,y)$ can be extended beyond  $y$ as segment, $y$ must be a branching point. By assumption branching points are $m$-negligible.  \hfill $\Box$
~\\

{\it Remark.} The example $(\R^n , \norm{.}_1,\lambda^n)$ shows that the $MCP(N,k)$-property is not strong enough to prevent a 'large' (with respect to $m$) amount of branching points, even if branching points  indicate infinite negative sectional curvature in Alexandrov sense. It is natural to ask which additional assumptions on $(X,d,m)$ inhibit a set of branching points with positive $m$-mass. For example,  $(X,d)$ admits no branching if it is a limit of Riemannian manifolds with uniform local lower sectional curvature bounds.\\

{\bf Acknowledgments.} Thanks to Jeff Cheeger for raising the question and to Karl-Theodor Sturm sending his preprints. Thanks also to Shin-ichi Ohta  in particular for bringing the preprint \cite{RM} to my attention. Ranjbar-Motlagh obtains  very similar results independent of the mass transportation approach but assumes a 'strong local doubling' property of $m$ instead. 

\def\cprime{$'$} \def\cprime{$'$} \def\cprime{$'$}

\end{document}